\documentclass[12pt,a4paper]{amsart}
\usepackage[top=3cm,bottom=3cm,left=2.8cm,right=2.8cm]{geometry}
\usepackage{cases}
\usepackage{graphicx}
\usepackage{stmaryrd}
\usepackage[all,pdf]{xy}
\usepackage{amssymb}
\usepackage{mathrsfs}
\usepackage{ulem}
\usepackage{xcolor}

\newtheorem{theorem}{Theorem}[section]

\newtheorem{lemma}[theorem]{Lemma}
\newtheorem{proposition}[theorem]{Proposition}
\newtheorem{definition}[theorem]{Definition}

\numberwithin{equation}{section}

\newtheorem{conjecture}{Conjecture}

       \newcommand{\beqm}{\begin{eqnarray*}}
       \newcommand{\eqm}{\end{eqnarray*}}

\newcommand{\fr}{\frac}

\begin{document}
\title{ Absolutely summing  Toeplitz operators on Fock spaces}

\author{Zhangjian Hu and Ermin Wang}

\address{Zhangjian Hu: Department of Mathematics, Huzhou Normal University, Huzhou 313000, Zhejiang  China}

\email{huzj@zjhu.edu.cn}

\address{Ermin Wang: School of Mathematics and Statistics, Lingnan Normal University, Zhanjiang 524048, Guangdong  China}

\email{wem0913@sina.com}

\thanks{}

\keywords{ Fock spaces; Absolutely Summing operators; Toeplitz operators}

\subjclass[2020]{30H20, 47B10}

\maketitle
\begin{abstract}
For $1\le p<\infty$, let $F^p_\varphi$ denote the Fock space on ${\mathbb C}^n$ with weight $\varphi$, where \(\varphi  \in  {\mathcal{C}}^{2}\left( {\mathbb{C}}^{n}\right)\) is  real-valued  and satisfies $m{\omega }_{0} \leq  d{d}^{c}\varphi \leq M{\omega }_{0}$ for two positive constants \(m\) and \(M\). Here \({\omega }_{0} = d{d}^{c}{\left| z\right| }^{2}\) is the Euclidean K\"{a}hler form on \({\mathbb{C}}^{n}\), and \({d}^{c} = \frac{\sqrt{-1}}{4}\left( {\bar{\partial } - \partial }\right)\).  In this paper, we completely characterize those positive Borel measures  $\mu$ on ${\mathbb C}^n$ for which the induced Toeplitz operators $T_\mu$ are $r$-summing on $F_{\varphi}^{p}$ for $r \ge 1$. 
\end{abstract}

 \section{  Introduction}
 The aim of this paper is to investigate $r$-summing property of Toeplitz operators on a class of Fock-type spaces for all $r\ge  1$.  In the following, let ${\mathbb C}^n$ be the $n$-dimensional complex Euclidean space. For two points $z=(z_1,\cdot\cdot\cdot,z_n)$ and $w=(w_1,\cdots,w_n)$ in ${\mathbb C}^n$, we write $\langle z,w\rangle=z_1\overline{w}_1+\cdots+z_n\overline{w}_n$
and $|z|=\sqrt{\langle z,z\rangle}.$ Write $B(z,\delta)$ for the ball centered at $z$ with radius $\delta>0$, that is, $B(z,\delta)=\{w\in {\mathbb C}^n: |w-z|<\delta\}$.
Throughout the paper, we always assume that \(\varphi  \in  {\mathcal{C}}^{2}\left( {\mathbb{C}}^{n}\right)\) is a real-valued function and satisfies
 $$
 m{\omega }_{0} \leq  d{d}^{c}\varphi  \leq M{\omega }_{0}
 $$
  for two positive constants \(m\) and \(M\) (write as $d{d}^{c}\varphi\simeq {\omega }_{0}$), where \({\omega }_{0} = d{d}^{c}{\left| z\right| }^{2}\) is the Euclidean K\"{a}hler form on \({\mathbb{C}}^{n}\), \({d}^{c} = \frac{\sqrt{-1}}{4}\left( {\bar{\partial } - \partial }\right)\).

Let  $dv$ be the Lebesgue volume measure on ${\mathbb C}^n$. For \(1\le  p < \infty\), we use $L^p$ to stand for the usual p-th Lebesgue space with the norm
$$
\|\cdot\|_{L^p}=\left(\int_{{\mathbb{C}}^{n}}{|\cdot| }^{p}dv \right)^{\frac{1}{p}} .
$$
The space \({L}^{p}_\varphi\) consists of all Lebesgue measurable functions \(f\) on $\mathbb{C}^n$ for which
\[
\|f\|_{L^p_\varphi } = {\left(\int_{{\mathbb{C}}^{n}}{\left| f(z) \right| }^{p}{e}^{-{p\varphi }(z) }dv(z) \right) }^{\frac{1}{p}} < \infty.
\]
Let \(H\left( {\mathbb{C}}^{n}\right)\) be the family of all holomorphic functions on \({\mathbb{C}}^{n}\). The  Fock space ${F}^{p}_\varphi$ is defined by
\[
{F}^{p}_\varphi = {L}^{p}_\varphi \cap H\left( {\mathbb{C}}^{n}\right)
\]
with inherited norm for \(1\le p < \infty\), and
$$
F_\varphi^\infty=\{f\in H({\mathbb{C}}^{n}): \| f\|_{F^\infty_\varphi}= \sup_{z\in {\mathbb{C}}^{n}}|f(z)|{e}^{-{\varphi }(z)}< \infty \}.
$$
 It is clear that \({F}^{p}_\varphi\) is a Banach space under the norm \(\|\cdot\|_{F^p_\varphi}\) for \(1 \leq  p \le \infty\).
The typical model for the weight function is $\varphi(z)=\frac \alpha 2|z|^2$, and the induced space is the classical Fock space denoted by $F^p_\alpha$, see \cite{Zh12} for more information.

Let \({K}_{\varphi }\left( {\cdot,\cdot}\right)\) be the Bergman kernel of ${F}^{2}_\varphi$, and set $k_{\varphi,z}(\cdot)=\frac{K_\varphi(\cdot,z)}{\sqrt{{K}_{\varphi }\left( {z,z}\right)}}$.
The orthogonal projection \(P_\varphi: {L}^2_\varphi  \rightarrow  {F}^{2}_\varphi\) can be represented as
\[
P_\varphi f(z) = \int_{\mathbb{C}^n}K_{\varphi}\left({z,w}\right) f(w) {e}^{-{2\varphi }\left( w\right) }dv(w).
\]
With this expression, $P_\varphi$ can be extended as  a bounded linear operator from ${L}^{p}_\varphi$ to ${F}^{p}_\varphi$ for $1\le p\le \infty$, and $P_\varphi f=f$ for $f\in {F}^{p}_\varphi$. For $1\le p<\infty$, the dual space of ${F}^{p}_\varphi$  can be identified with ${F}^{p'}_\varphi$ under the pairing $\langle f, g \rangle_\varphi$, which is defined as
$$
\langle f, g \rangle_\varphi=\int_{{\mathbb{C}}^{n}}f(z)\overline{g(z)} {e}^{-{2\varphi }(z) }dv(z),
$$
where $p'$ is the conjugate exponent of $p$, that is, $\frac 1p+ \frac 1{p'} =1$. See \cite{SV12} for details.

Given some positive Borel measure $\mu$ on $\mathbb{C}^n$ (denoted by $\mu \ge 0$), the Toeplitz operator $T_\mu$ with symbol $\mu$ is defined by
\[
T_\mu f(z) = \int_{\mathbb{C}^n}K_{\varphi}\left({z,w}\right) f(w) {e}^{-{2\varphi }\left( w\right) }d\mu(w).
\]
And also, we set $$\widetilde{\mu}(z)=\int_{{\mathbb C}^n} \left| k_{\varphi, z}(w)\right|^2 d\mu(w)$$ to be the Berezin transform of $\mu$, and $$\widehat{\mu}_\delta(z)= \mu(B(z, \delta))=\int_{B(z, \delta)}d\mu(w)$$ with  $\delta>0$ fixed.

\begin{definition} \label{def1}
Suppose \(1 \leq  r< \infty\) , a linear operator \(T : X \rightarrow  Y\) between Banach spaces \(X\) and \(Y\) is said to be \(r\)-summing if there exists \(C \geq  0\) such that
\begin{equation}\label{Main-def}
{\left( \mathop{\sum }\limits_{{k = 1}}^{n}{\begin{Vmatrix}T{x}_{k}\end{Vmatrix}}_{Y}^{r}\right) }^{\fr 1r} \leq  C\mathop{\sup }\limits_{{{x}^{ * } \in  {B}_{{X}^{ * }}}}{\left( \mathop{\sum }\limits_{{k = 1}}^{n}{\left| {x}^{*}({x}_{k}) \right| }^{r}\right) }^{\fr 1r} 
\end{equation}
for every finite sequence \({\left\{  {x}_{k}\right\}  }_{1 \leq  k \leq  n} \subset  X\), where ${X}^{ * }$ denotes the dual space of $X$, and ${B}_{{X}^{*}}$ denotes the closed unit ball $\{x^*\in {X}^{*}: \|x^*\|\le 1\}$ of  ${X}^{ * }$.
\end{definition}

We shall write $\Pi_r(X,Y)$ for the set of all $r$-summing operators from $X$ to $Y$.
The \(r\)-summing norm of \(T\), denoted by \({\pi }_{r}\left( T: X\to Y\right)\), is the least possible constant \(C\) in (\ref{Main-def}). In the following, if no  confusion occurs, we will abbreviate $\pi_r(T: X\to Y)$ as $\pi_r(T)$.\\

 The study of absolutely summing operators traces its origins to Grothendieck's pioneering contributions in the 1950s. In his work on nuclear spaces, he introduced the basic properties of these operators (see \cite{Gr53, Gr55}). A major advancement in the theory came from Pietsch, who explicitly defined the class of $r$-absolutely summing operators for all $1 \leq r < \infty$ in \cite{Pi67}. This development significantly generalised Grothendieck's original concept and led to the establishment of many foundational properties of such operators. We refer to \cite{Pi72, Pi78} for more information. More recently, considerable attention has been directed towards the study of absolutely summing Carleson embeddings on various analytic function spaces. For example, \cite{LR18} explores the Hardy spaces $H^p$, \cite{HJL24} focuses on the classical Bergman spaces $A^p(\mathbb{D})$, and \cite{CHW24} investigates weighted Fock spaces $F^p_{\alpha, \omega}$ with $A_\infty$-type weights. Furthermore, absolutely summing weighted composition operators on Bloch spaces, as well as Volterra operators on Bergman and Bloch spaces, have been explored in \cite{FL22, JL24}. However, it is important to note that all of the aforementioned research is confined to function spaces of a single variable.

Toeplitz operators are widely recognized for their significant applications in fields such as signal processing and quantum theory. They also play a crucial role in the 
operator theory of holomorphic function spaces. These operators have been extensively studied on classical function spaces like Hardy and Bergman spaces. In the context of Fock spaces, the systematic study of Toeplitz operators was initiated by Berger and Coburn \cite{BC86, BC87}. Later, Isralowitz and Zhu \cite{IZ10} established criteria for the boundedness, compactness, and Schatten class membership of Toeplitz operators with positive symbols on the classical Fock space $F^2_{\alpha}$. Subsequent research extended these studies to various generalizations of Fock spaces (see \cite{HL11, HL14, HL17, IVW15, OP16} and the references therein).

It is worth noting that the existing research on Toeplitz operators has mainly focused on their boundedness, compactness, Schatten class properties, and related algebraic structures. The theory of absolutely summing Toeplitz operators on holomorphic function spaces--including, but not limited to, Fock spaces--remains largely unexplored.

The present work aims to fill this gap by providing a comprehensive characterization of positive Toeplitz operators that are $r$-summing on Fock spaces $F^p_\varphi$. This contribution not only enriches the general theory of operator ideals in analytic function spaces but also has potential implications in areas where Toeplitz operators serve as mathematical models, such as time-frequency analysis, quantization in mathematical physics, and complex geometry. From a methodological perspective, our results reveal new connections between the Berezin transform, Carleson-type conditions, and the geometric structure of Fock spaces, thereby offering a framework for further developments in related settings.

We need some more  notation to state  our main result.  Two quantities $A$ and $B$ are equivalent, denoted by  $A\simeq B$, if  $C^{-1}A \le B \le CA$ for some $C>0$. The conclusion in the following Theorem \ref{THM2} is new even in the simplest case that $n=1$ and $\varphi(z)= \frac \alpha 2 |z|^2$ with $\alpha>0$ fixed.

\begin{theorem}\label{THM2}
 Let \(\varphi\in C^2(\mathbb{C}^n)\) satisfying $dd^c\varphi\simeq \omega_0$, and let $\mu$ be a positive Borel measure on \(\mathbb{C}^n\).

{\rm (1)} If \(1\le p \le 2\), then for  \(r \geq  1\), \(T_\mu: F_\varphi^p \rightarrow  F_\varphi^p\) is \(r\)-summing if and only if \(\widetilde{\mu} \in  {L}^{2}\).

{\rm (2)} If \(2\le p <\infty\), then

 \ \ \, {\rm (i)} for $1 \leq  r \leq p'$, \({T}_{\mu} : {F}_{\varphi}^{p} \rightarrow  {F}_{\varphi }^{p}\) is \(r\)-summing if and only if \(\widetilde{\mu} \in  {L}^{p'}\).

 \ \ \  {\rm (ii)} for ${p'} \le  r \leq  p$, \({T}_{\mu} : {F}_{\varphi}^{p} \rightarrow  {F}_{\varphi }^{p}\) is \(r\)-summing if and only if \(\widetilde{\mu} \in  {L}^{r}\).

 \  \,  {\rm (iii)} for $p\le r<\infty$, \({T}_{\mu} : {F}_{\varphi}^{p} \rightarrow  {F}_{\varphi }^{p}\) is \(r\)-summing if and only if \(\widetilde{\mu} \in  {L}^{p}\).\\
Moreover, we have the quantity equivalence
\begin{equation}\label{main-thm} \pi_r(T_\mu : F_{\varphi}^{p} \rightarrow  F_{\varphi}^{p}) \simeq
\left\{
  \begin{array}{ll}
   \|\widetilde{\mu}\|_{L^{2}}, \ \ \ \ \textrm{if} \ 1\le p\le 2; \\
    \|\widetilde{\mu}\|_{L^{p'}}, \ \ \ \textrm{if} \ p\ge 2, \ \textrm{ and }\ 1\le r\le p';\\
    \|\widetilde{\mu}\|_{L^{r}}, \, \ \ \ \textrm{if} \ p\ge 2, \ \textrm{ and }\ p'\le r\le p;\\
    \|\widetilde{\mu}\|_{L^{p}}, \, \ \ \ \textrm{if} \ p\ge 2, \ \textrm{ and }\ p\le r<\infty.
  \end{array}
\right.
\end{equation}
And all statements above remain true if $\widetilde{\mu}$ is replaced by $\widehat{\mu}_\delta$ with $\delta>0$ fixed.
\end{theorem}

The remainder of this paper is organized as follows. In Section 2, we introduce preliminary results, establish notation, and review several standard facts that will be referenced throughout the paper. Section 3 is dedicated to the proof of our main theorem. Our approach to characterizing the $r$-summing Toeplitz operators departs significantly from the methods used in \cite{CHW24, FL22, HJL24, JL24, LR18}, which focus on absolutely summing Carleson embeddings, composition operators, and Volterra operators. To analyze the $r$-summability of the operator $T_\mu: F^p_\varphi \to F^p_\varphi$ for $1 \le p \le 2$, we employ the interpolation technique based on tensor products of Banach spaces. For the case $p \ge 2$, we make use of the connections between absolutely summing operators and order boundedness, combining them with Grothendieck's theorem and Pietsch's domination theorem. This methodology enables us to examine the absolutely summing property of $T_\mu$ even when the exponent of the underlying holomorphic function space is 1--a case not previously addressed in the literature. In Section 4, we discuss potential extensions of our work, including the study of absolutely summing Toeplitz operators between different Fock spaces. Finally, we propose a conjecture regarding $r$-summing Toeplitz operators on Bergman spaces.

 Throughout this paper, we use $C$ to denote a positive constant whose value may change from line to line, but which does not depend on the functions under consideration. For two quantities $A$ and $B$, we write  $A\lesssim B$ if there exists a constant $C$ such that   $A\leq CB$. Again, $A$ and $B$ are said to be equivalent, denoted  $A\simeq B$, if  $C^{-1}A \le B \le CA$ for some $C>0$.

 \section{Preliminaries}
\subsection{Some notations and  estimates.}
The following lemma  provides some estimates for the Bergman kernel of $F^2_\varphi$, which can be found in  \cite{SV12}.
\begin{lemma}\label{lem1}
Let $K_\varphi(\cdot,\cdot) $ be the Bergman kernel of $F^2_\varphi$, then we have the following estimates:

{\rm (1)} There exists $\theta>0$ such that
 \begin{equation}\label{formula1}
|K_\varphi(z,u) |e^{-{\varphi(z)-\varphi(u)}}\lesssim e^{-{\theta|z-u|}}, \  z,u\in \mathbb{C}^n.
\end{equation}

{\rm (2)}  There exists some $r>0$ such that
$$
|K_\varphi(z,u) |e^{-{\varphi(z)-\varphi(u)}}\gtrsim  1 \  {\rm whenever} \ z\in \mathbb{C}^n, u\in B(z,r).
$$

{\rm (3)} For $0<p\le \infty$,
$$
\|K_\varphi(\cdot,z) \|_{L^p_\varphi}\simeq e^{\varphi(z)}\simeq \sqrt{K_\varphi(z,z)}, \  z\in \mathbb{C}^n.
$$
\end{lemma}

 Given some sequence $\{a_k\}_{k=1}^{\infty}$ in $\mathbb{C}^n$ and $\delta>0$,  we call $\{a_k\}_{k=1}^{\infty}$  a $\delta$-lattice if the balls  $\left\{B(a_k,\delta)\right\}_{k=1}^{\infty}$ cover $\mathbb{C}^n$ and  $\left\{B\left(a_k,\frac{\delta}{2\sqrt{n}}\right)\right\}_{k=1}^{\infty}$ are pairwise  disjoint. A typical model of a $\delta$-lattice is the sequence
 $$
 \left\{\frac{\delta}{\sqrt{n}}(m_1+ {\rm i}k_1,m_2+ {\rm i}k_2, ..., m_n+ {\rm i}k_n)\in \mathbb{C}^n: m_j, k_j \in \mathbb{Z}, j=1,2,...,n \right\}.
 $$
With the above hypothesis,  there exists an integer $N$ so that for any $\delta$-lattice $\{a_k\}_{k=1}^{\infty}$,
$$
   1\le \sum_{k=1}^\infty \chi_{B(a_k,2\delta)}(z)\le N
$$
for $ z\in {\mathbb{C}^n}$, where $\chi_{E}$ is the characteristic function of a subset $E$ of $\mathbb{C}^n$.

The following lemma comes from \cite{HL14}.
\begin{lemma}\label{lem3}
Let $\{a_k\}_{k=1}^\infty$ be a $\delta$-lattice. Then for  $0<p\le \infty$ and $\{c_k\}\in l^p$, the function  $f$ defined by
$$
f(z)=\sum^{\infty}_{k=1}c_kk_{\varphi,a_k}(z), \ z\in \mathbb{C}^n
$$
belongs to $F^p_\varphi$ with $\|f \|_{L^p_\varphi}\lesssim \|\{c_k\} \|_{l^p}$.
\end{lemma}

 For $0<p,q<\infty$, a positive measure $\mu$ is called a $(p,q)$-Carleson measure if the embedding $${\rm Id}: {F}^{p}_\varphi \to {L}^{q}_\varphi(\mu)$$ is bounded.  In \cite{HL14}, Hu and Lv provided a complete characterization of such measures, which consequently resolved the question of boundedness of the Toeplitz operator  $T_\mu $ from ${F}^{p}_\varphi$ to $F^{q}_\varphi$ for all $0<p,q<\infty$. Below, we summarize the relevant results required for our study.

\begin{lemma}\label{lem2.3}
Suppose $\mu\ge 0$. Then we have the following statements:

{\rm (1)} For $0<p\le \infty$,  $\widetilde{\mu}\in L^p\Leftrightarrow\widehat{\mu}_\delta\in L^p\Leftrightarrow\{\widehat{\mu}_\delta(a_k)\}\in l^p$ for any $\delta$-lattice $\{a_k\}$. Moreover,
$$
\|\widetilde{\mu}\|_{L^p}\simeq\|\widehat{\mu}_\delta\|_{L^p}\simeq\|\{\widehat{\mu}_\delta(a_k)\}\|_{l^p}.
$$

{\rm (2)}  For $0<q<p<\infty$, $T_\mu: {F}^{p}_\varphi\to F^{q}_\varphi$ is bounded if and only if $\widetilde{\mu}\in L^{\frac{pq}{p-q}}$. Moreover,
$$
\|T_\mu\|_{F^p_\varphi\to F^q_\varphi}\simeq\|\widetilde{\mu}\|_{L^{\frac{pq}{p-q}}}.
$$
\end{lemma}

\subsection{Absolutely summing operators.}

For a comprehensive treatment of $r$-summing operators, we refer readers to \cite{DJT95, Pe77, Wo91}. In what follows, we recall several basic properties of these operators that are essential to our study, all of which can be found in \cite{DJT95}.

{\rm (1)} \textbf{Ideal property of $p$-summing operators} {\rm[p.37]}: For any $r \geq  1$, the class ${\Pi }_{r}\left( {X,Y}\right)$ forms an operator ideal between Banach spaces: for any $T \in  {\Pi }_{r}\left( {X,Y}\right)$ , and for any two Banach spaces ${X}_{0},{Y}_{0}$ such that both $S : {X}_{0} \rightarrow  X$ and $U : Y \rightarrow  {Y}_{0}$ are bounded linear operators, we have ${UTS} \in  {\Pi }_{r}\left( {{X}_{0},{Y}_{0}}\right)$ with
\begin{equation}\label{ideal-a}
{\pi }_{r}\left( {UTS}\right)  \leq  \|U\| \cdot{\pi }_{r}( T)\cdot \| S\|,
\end{equation}

{\rm (2)} \textbf{Inclusion relations} {\rm[p.39]}: If  $1 \le p<q<\infty$, then ${\Pi }_{p}\left( {X,Y}\right)\subset{\Pi }_{q}\left( {X,Y}\right)$. Moreover, for $T \in  {\Pi }_{p}\left( {X,Y}\right)$, we have
\begin{equation}\label{inclusion}
{\pi }_{q}(T)\le {\pi }_{p}(T) \ \ \textrm{for } \ \ q\ge p.
\end{equation}

{\rm (3)} \textbf{Cotype property}:
 Let $\{r_k\}_{k=1}^\infty$  be  the sequence of Rademacher functions on  $[0, 1]$, defined by
$$ r_{1}(t)=
\left\{
  \begin{array}{ll}
    1, \ \ \ \ \ \textrm{if} \ 0\leq t-[t]<\frac{1}{2}; \\
    -1, \ \ \ \textrm{if} \ \frac{1}{2}\leq t-[t]<1,
  \end{array}
\right.
$$
and for $k=1,2,\dots,$
\begin{equation}\label{rademacher}
r_{k+1}(t)=r_{1}(2^{k}t).
 \end{equation}
We will use a powerful inequality about the  Rademacher functions in this paper. Precisely, for $0<l<\infty$,  there exist positive constants $C_{1}$  and $C_{2}$ depending only on $l$ such that
\begin{equation}\label{Khinchine-1}
C_{1}
\left(\sum_{k=1}^{m}|b_{k}|^{2}\right)^{\frac{l}{2}}\leq\int_{0}^{1}\left|\sum_{k=1}^{m}b_{k}r_{k}(t)\right|^{l}dt \leq C_{2}
\left(\sum_{k=1}^{m}|b_{k}|^{2}\right)^{\frac{l}{2}}
\end{equation}
for all $m\geq 1$ and complex numbers $b_{1}, b_{2}, \ldots, b_{m}$. The above inequality is known as the Khintchine inequality. See \cite{DJT95} for the details.

 We say that  a Banach space $X$ has cotype $q$ if there is a constant $\kappa > 0$ such that no matter how we select finitely many vectors $x_1, x_2, \cdots, x_m$  from $X$,
$$
    \left(\sum_{j=1}^m \|x_j\|^q\right)^{\frac 1 q}\le \kappa \left( \int_0^1 \left\| \sum_{j=1}^m r_j(t) x_j\right\|^2dt\right)^{\frac 12},
$$
where $r_j$ is as in (\ref{rademacher}).

It is well known that,  for $1\le p<\infty$, the Lebesgue  space
$L^p(\Omega,\mu)$ has cotype ${\rm max}\{p,2\}$,
see  \cite[p.219]{DJT95}.   For Banach spaces  $X$ and $Y$,
if $X$ and $Y$ both have cotype 2,   \cite[Corollary 11.16]{DJT95} tells us that
\begin{equation} \label{cotype}
{\Pi }_{r}\left( {X,Y}\right)={\Pi}_{1}\left( {X,Y}\right)
\end{equation}
with the equivalence $\pi_r(T)\simeq\pi_1(T)$ for $1\le r<\infty$.
We shall use the following well known fact several times about $r$-summing operators, which gives an equivalent description of $r$-summing operators, see \cite{DJT95} or \cite{LR18}.
\begin{lemma}\label{lem8}
  Let \(X,Y\) be Banach spaces and \(r \geq  1\). Then a bounded operator $T: X\to Y $ is $r$-summing if and only if there is a constant $C>0$ such that  for every $X$-valued random variable $F$ on any measure space $(\Omega, \nu)$, we have
\[
\left(\int_{\Omega }\|T \circ  F\|_{Y}^{r}{d\nu}\right)^{\fr 1r} \leq C\mathop{\sup }\limits_{{\xi  \in  {B}_{{X}^{ * }}}}\left(\int_{\Omega }{\left| \xi  \circ  F\right| }^{r}{d\nu}\right)^{\fr 1r}.
\]
Moreover, the best constant $C$ is  ${\pi }_{r}{\left( T\right) }$.
\end{lemma}

Order bounded operators are closely related to $r$-summing operators and play an equally important role in our study.   Recall that for $p\ge 1$, a Banach space operator $T: X\to {L}^{p}(\mu)$ is order bounded if $T(B_X)$ is an order bounded subset of $L^p(\mu)$, that is, there exists a non-negative function $h\in L^p(\mu)$ such that for each $f\in B_X$, $|Tf|\le h$ $\mu$-almost everywhere. The following lemma, which follows from Propositions 5.5 and 5.18 in \cite{DJT95}, establishes a connection between these two classes of operators. 

\begin{lemma}\label{lem66}
Let $X$ be a Banach space, $p \geq  1$, and let $(\Omega ,\Sigma ,m)$ be a measure space. If the operator $T : X \rightarrow  {L}^{p}\left( {\Omega ,m}\right)$ is order bounded, then it is p-summing with
$$
{\pi }_{p}\left( T\right)  \leq  {\left\|\mathop{\sup }\limits_{f \in  B_X}|Tf| \right\|}_{{L}^{p}\left( {\Omega ,m}\right) }.
$$
\end{lemma}

\subsection{Injective tensor product.}  In our analysis we will use the  interpolation technique on tensor product of Banach spaces.  Here are some  necessary notations and results which can be consulted from \cite{Ry02}.

 Let $X$ be a Banach space. For $p\ge 1$, the vector space $l^p(X)$ consists of all sequences $\{x_j\}_{j=1}^\infty$ in $X$ such that
$$
\left \|\{x_j\}\right \|_{l^p}=\left(\sum^{\infty}_{j=1}\|x_j\|^p\right)^{\fr 1 p}< \infty.
$$
We denote by $l^{p,w}(X)$ the linear space of the sequences $\{x_j\}_{j=1}^\infty$ in $X$ such that
$$
\|\{x_j\}\|_{w,p}=\sup_{x^\ast\in B_{X^
{\ast}}}\|\{x^\ast(x_j)\}\|_{l^p}<\infty.
$$
As usual, $l^{p,u}(X)$ denotes the subspace of $l^{p,w}(X)$ consisting of those sequences $\{x_j\}\in l^{p,w}(X)$ such that
$$
\lim_{k\to \infty}\|\{x_{j+k}\}_{j=1}^\infty\|_{w,p}=0.
$$
Both $l^{p,w}(X)$ and $l^{p,u}(X)$ are Banach spaces.

Given  Banach spaces $X$ and $Y$, let $\mathcal{B}(X, Y)$ be the collection of all bounded linear operators from $X$ to $Y$. For $T\in {\mathcal B}(X, Y)$, the correspondence
$$
\widehat{T}: \{x_n\}_n\mapsto \{Tx_n\}_n
$$
always induces a bounded linear operator, that is, $\widehat{T} \in {\mathcal B} \left(l^{p,w}(X),  l^{p,u}(Y)\right)$.   It follows from \cite[Proposition 2.1]{DJT95} that
$$
T\in \Pi_r(X, Y)  \ \textrm{ if and only if } \ \ \widehat{T}\in \mathcal{B}(l^{r,w}(X), l^r(Y)).
$$
 With the same approach as in \cite{DJT95} we see that
\begin{equation}\label{w-to-str}
T\in \Pi_r(X, Y) \ \textrm{ if and only if } \ \  \widehat{T}\in \mathcal{B}(l^{r,u}(X), l^r(Y))
\end{equation}
with the norm equivalence
\begin{equation}\label{w-to-str-a}
\pi_r(T)\simeq\|\widehat{T}\|_{l^{r,u}(X)\to l^r(Y)}.
\end{equation}

 For two vector spaces $X$ and $Y$, let $X \otimes Y$  be the  tensor product as  in \cite{Ry02}. A typical tensor $u$ in $X \otimes Y$ has the form
\begin{equation}\label{tensor-1}
      u= \sum_{j=1}^m x_j \otimes y_j,  \ x_j\in X,  \ y_j\in Y.
\end{equation}
  When $X$ and $Y$ are Banach spaces, we define the injective norm on  $X\otimes Y$ as
 $$
    \varepsilon(u)= \sup \left \{ \left|\sum_{j=1}^n x^*(x_j) y^*(y_j)\right|: x^*\in B_{X^\ast}, y^*\in B_{Y^\ast} \right \},
 $$
 where $u$ is as (\ref{tensor-1}).
The completion of $X \otimes Y$ endowed with this norm is called the injective tensor product of $X$ and $Y$,  denoted by $X\widetilde{\otimes}_\varepsilon Y$. See \cite{DC00} for the details.

Now we consider the product $l^p \widetilde{\otimes}_\varepsilon X$ with $1\le p<\infty$. Define the map $L$ from $l^p  \otimes X$ to $l^{p,w}(X)$ as
$$
    Lu= \left\{ \sum_{j=1}^m a_{j, k} x_j \right \}_{k=1}^\infty,  u = \sum_{j=1}^m a_j\otimes x_j\in l^p  \otimes X\,  \textrm{with} \, a_j= \{a_{j, k}\}_{k=1}^\infty.
$$
With the same approach as that in \cite[Examples 3.4]{Ry02} we know that $L$ is an isometric isomorphism between $l^p \widetilde{\otimes}_\varepsilon X$ and
$l^{p,u} (X)$. In other words, as mentioned in  \cite[page 3988]{JMP08}, we  have the natural isometric identification 
\begin{equation}\label{tensor-2}
    l^p \widetilde{\otimes}_\varepsilon X= l^{p,u} (X)  \ \ \textrm{ for } \ \ 1\le p<\infty.
\end{equation}

\subsection{Complex interpolation.}
We are going to prove the main result for $1\le p\le 2$ by using interpolation techniques. We assume the reader is familiar with the complex interpolation method of Calder\'{o}n (see \cite{BL78}).

If $X_0$ and $X_1$ are compatible  Banach spaces and $0<\theta<1$, we have the interpolation functor $C_\theta$ so that $ C_\theta(X_0, X_1)$ (denoted by $[X_0, X_1]_\theta$)  is an interpolation space between $X_0$ and $X_1$. Furthermore, if $X_0, X_1$ and $Y_0, Y_1$ are two compatible pairs of Banach spaces  and
$$
    T: X_0+X_1\to Y_0+Y_1
$$
is a  linear mapping such that $T$ maps $X_j$ boundedly into $Y_j$  ($j=0,1$), then
$$
T:[X_0, X_1]_\theta\to [Y_0, Y_1]_\theta \ \  \textrm {boundedly for all } \ \ \theta\in (0, 1)
$$
with the norm estimate
\begin{equation}\label{interpolation-a}
\|T\|_{[X_0, X_1]_\theta \to [Y_0, Y_1]_\theta}\le M_0^{1-\theta} M_1^\theta,
\end{equation}
where
$$
   M_0= \|T\|_{X_0\to Y_0} \ \textrm{ and }\  M_1= \|T\|_{X_1\to Y_1}.
$$
See \cite[Theorem 4.1.2]{BL78} and \cite[Theorem 2.4]{Zh07}.

We will apply the complex interpolation on Banach lattices. For the basic theory of Banach lattices we refer the reader to \cite{LT79}. As in \cite{Ko91}, an interpolation pair of Banach spaces $(X_0, X_1)$ is called an interpolation pair of Banach lattices if
 $X_0$ and $X_1$ are Banach lattices which can be embedded as lattice-ideals in some  $L^0(\Omega, \Sigma, \mu)$, the space of all measurable functions on a $\sigma$-finite measure space $(\Omega, \Sigma, \mu)$.  Moreover, the set of simple functions
 $$
    S= \left\{\sum_{j=1}^m \alpha_j \chi_{\Omega_j}: m \in \mathbb N, \alpha_j\in \mathbb C, \Omega_j\in \Sigma, \mu(\Omega_j)<\infty\right\}
 $$
 is dense in both $X_0$ and $X_1$.

  A Banach lattice $X$ is called 2-concave if there exists some $M>0$ such that
$$
\left(\sum^{m}_{j=1}\|x_j\|^2_X\right)^{\fr 1 2}\le M\left\|\left(\sum^{m}_{j=1}|x_j|^2\right)^{\fr 1 2}\right\|_{X}
$$
for every choice of vectors $\{x_j\}_{j=1}^m$ in $X$.

The following proposition  is copied from \cite[Theorem 4.2]{Ko91}.

 \begin{proposition}\label{Ko91-1} Let $(X_0, X_1)$ and $(Y_0, Y_1)$ be two interpolation pairs, and $X_0^{\ast}, X_1^{\ast}$ are type 2 spaces, $Y_0, Y_1$ are 2-concave Banach lattices. Then, $(X_0 \widetilde{\otimes}_\varepsilon Y_0, X_1\widetilde{\otimes}_\varepsilon Y_1)$ is an interpolation pair and for $0<\theta<1$ we have
 $$
     [X_0 \widetilde{\otimes}_\varepsilon Y_0, X_1\widetilde{\otimes}_\varepsilon Y_1]_\theta = [X_0, X_1]_\theta   \widetilde{\otimes}_\varepsilon [ Y_0, Y_1]_\theta.
 $$
 \end{proposition}

\section{Proof of Theorem 1.2}

In this section, we present the proof of Theorem \ref{THM2},  which is divided into four cases.

\subsection{Case 1: $1\le p \le2$.}

We prove the sufficiency first. Suppose \(\widetilde{\mu} \in L^ {2}\). We claim that
both \({T}_{\mu} : {F}_{\varphi}^{1} \rightarrow  {F}_{\varphi }^{1}\) and \({T}_{\mu} : {F}_{\varphi}^{2} \rightarrow  {F}_{\varphi }^{2}\) are \(1\)-summing, and
\begin{equation}\label{proof-a}
\pi_1(T_\mu: {F}_{\varphi}^{1} \rightarrow  {F}_{\varphi }^{1})\lesssim \|\widetilde{\mu}\|_{L^ {2}}, \ \ \pi_1(T_\mu: {F}_{\varphi}^{2} \rightarrow  {F}_{\varphi }^{2})\lesssim \|\widetilde{\mu}\|_{L^ {2}}.
\end{equation}
In fact, with the norm estimate $\|f\|_{F^2_\alpha}\lesssim \|f\|_{F^1_\alpha}$  we know
\beqm
 \mathrm{Id}_1: {F}_{\varphi}^{1} &\to & {F}_{\varphi }^{2} \\
f &\mapsto &f
\eqm
is bounded. Consider the chain of maps:
$$
L^1_\varphi\stackrel{P_\varphi}{\longrightarrow} F^1_\varphi\stackrel{{\rm Id}_1}{\longrightarrow} F^2_\varphi\stackrel{{\rm Id}_2}{\longrightarrow}L^2_\varphi.
$$
It is obvious that $J={\rm Id}_2\circ {\rm Id}_1\circ P_\varphi$ is bounded from $L_{\varphi}^{1}$ to ${L}_{\varphi }^{2}$ with $\|J\|_{L_\varphi^1 \rightarrow  {L}_{\varphi }^{2}}\lesssim 1$. By \cite[Theorem 3.4]{DJT95} we get $J: {L}_{\varphi}^{1} \rightarrow  {L}_{\varphi }^{2}$ is 1-summing with
$$\pi_1(J: L_\varphi^1 \rightarrow  L_\varphi ^2)\lesssim \|J\|_{ L_\varphi^1 \rightarrow  L_\varphi ^2}.$$
Then, by \cite[p.37]{DJT95} we get the restriction map $J|_{{F}_{\varphi}^{1}}: F_{\varphi}^{1} \rightarrow  {L}_{\varphi }^{2}$ is also 1-summing with
$$
\pi_1(J|_{{F}_{\varphi}^{1}}: F_\varphi^1 \rightarrow  L_\varphi ^2)\le \pi_1(J: L_\varphi^1 \rightarrow  L_\varphi ^2)\lesssim 1.
$$
Hence, (\ref{ideal-a}) shows
$$
{\rm Id}_1(F_\varphi^1 \rightarrow  F_\varphi ^2)=P_\varphi({L}_{\varphi}^{2} \rightarrow  {F}_{\varphi }^{2})\circ J|_{{F}_{\varphi}^{1}}({F}_{\varphi}^{1} \rightarrow  {L}_{\varphi }^{2})
$$
 is 1-summing, or
 \begin{equation}\label{1-to-2}
\pi_1({\rm Id}_1: F_\varphi^1 \rightarrow  F_\varphi ^2)< \infty.
 \end{equation}
 Meanwhile, Lemma \ref{lem2.3} (2) tells us  ${T}_{\mu} : {F}_{\varphi}^{2} \rightarrow  {F}_{\varphi }^{1}$ is bounded with $\|T_\mu\|_{F_\varphi^2 \rightarrow  {F}_{\varphi }^{1}}\lesssim \|\widetilde{\mu}\|_{L^ {2}}$.
 Therefore, (\ref{ideal-a}) turns out
$$
{T}_{\mu}({F}_{\varphi}^{1} \rightarrow  {F}_{\varphi }^{1})=T_\mu(F_\varphi^2 \rightarrow  {F}_{\varphi }^{1})\circ \mathrm{Id}_1 ({F}_{\varphi}^{1} \rightarrow  {F}_{\varphi }^{2})
$$
is \(1\)-summing, and
\begin{equation}\label{proof-b}
\pi_1({T}_{\mu}:{F}_{\varphi}^{1} \rightarrow  {F}_{\varphi }^{1})\le \|T_\mu\|_{F_\varphi^2 \rightarrow  {F}_{\varphi }^{1}}\cdot\pi_1(\mathrm{Id}_1: {F}_{\varphi}^{1} \rightarrow  {F}_{\varphi }^{2})\lesssim\|\widetilde{\mu}\|_{L^ {2}}.
\end{equation}
Similarly,
$$
{T}_{\mu}({F}_{\varphi}^{2} \rightarrow  {F}_{\varphi }^{2})=\mathrm{Id}_1({F}_{\varphi}^{1} \rightarrow  {F}_{\varphi }^{2})\circ T_\mu(F_\varphi^2 \rightarrow  {F}_{\varphi }^{1})
$$
is \(1\)-summing with
\begin{equation}\label{proof-c}
\pi_1({T}_{\mu}:{F}_{\varphi}^{2} \rightarrow  {F}_{\varphi }^{2})\lesssim \|\widetilde{\mu}\|_{L^ {2}}.
\end{equation}
 The estimates
(\ref{proof-a}) follow from  (\ref{proof-b}) and (\ref{proof-c}).

Now we show, for \(1< p <  2\),  \({T}_{\mu} : {F}_{\varphi}^{p} \rightarrow  {F}_{\varphi }^{p}\) is \(r\)-summing  with the estimate
\begin{equation}\label{proof-d}
\pi_r(T_\mu)\lesssim \|\widetilde{\mu}\|_{L^2}.
\end{equation}
Since the cotype of \({F}_{\varphi}^{p}\) is 2 for \(1< p<  2\), by (\ref{cotype})  it suffices to prove  $T_{\mu}\in {\Pi_2} ({F}_{\varphi}^{p}, {F}_{\varphi }^{p})$ for \(1\le p \le  2\). We are going to use interpolation to get the desired conclusion since we already know $T_{\mu}\in {\Pi_2} ({F}_{\varphi}^{p}, {F}_{\varphi }^{p})$ for $p=1, 2$.

Consider the chain of maps
$$
L^p_\varphi\stackrel{P_\varphi}{\longrightarrow} F^p_\varphi\stackrel{T_\mu}{\longrightarrow} F^p_\varphi\stackrel{{\rm Id}_3}{\longrightarrow}L^p_\varphi,
$$
 and set $T ={\rm Id}_3\circ T_\mu\circ P_\varphi: L_{\varphi }^{p}\to L_{\varphi }^{p}$.
We claim that $T_{\mu}\in {\Pi_r} ({F}_{\varphi}^{p}, {F}_{\varphi }^{p})$ if and only if $T\in {\Pi_r} ({L}_{\varphi}^{p}, {L}_{\varphi }^{p})$ and
\begin{equation}\label{342}
\pi_r(T_\mu)  \simeq  \pi_r(T).
\end{equation}
In fact, if $T\in {\Pi_r} ({L}_{\varphi}^{p}, {L}_{\varphi }^{p})$, then the restriction map $T|_{F_{\varphi }^p}: F_{\varphi }^p\to L_{\varphi }^p$ is also $r$-summing with $\pi_r(T|_{F_{\varphi }^p})\le \pi_r(T).$ Thus we get
$$
T_\mu=P_\varphi\circ T|_{F_\varphi^p}
$$
is $r$-summing with
$$
\pi_r(T_\mu)\le \|P_\varphi\|_{L_\varphi^p \rightarrow F_\varphi ^p}\cdot\pi_r(T|_{F_\varphi ^p})\lesssim \pi_r(T).
$$
 Conversely, if $T_{\mu}\in {\Pi_r} ({F}_{\varphi}^{p}, {F}_{\varphi }^{p})$, since both ${\rm Id}_3: F_{\varphi }^{p}\to L_{\varphi }^{p}$ and $P_\varphi: L_\varphi^p \rightarrow F_\varphi ^p$ are bounded, the ideal property of $r$-summing operators implies that $T={\rm Id}_3\circ T_\mu\circ P_\varphi$ is $r$-summing with
$$
\pi_r(T)\le \|{\rm Id}_3\|_{F_\varphi^p \rightarrow L_\varphi ^p}\cdot\pi_r(T_\mu)\cdot\|P_\varphi\|_{L_\varphi^p \rightarrow F_\varphi ^p}\lesssim \pi_r(T_\mu).
$$
 Thus (\ref{342}) holds. Meanwhile, (\ref{w-to-str}) tells us $
T\in \Pi_2(L_{\varphi }^{p}, L_{\varphi }^{p})$ if and only if $\widehat{T}\in \mathcal{B}\left(l^{2,u}(L_{\varphi }^{p}), l^2(L_{\varphi }^{p})\right).
$
Moreover,
\begin{equation}\label{352}
\|\widehat{T}\|_{ l^{2,u}(L_{\varphi }^{p})\to l^2(L_{\varphi }^{p})}\simeq\pi_2(T: L_{\varphi }^{p}\to L_{\varphi }^{p}).
\end{equation}
This and (\ref{proof-a}) imply
$$
\|\widehat{T}\|_{ l^{2,u}(L_{\varphi }^{1})\to l^2(L_{\varphi }^{1})}<\infty\ \ \textrm{and} \ \ \|\widehat{T}\|_{ l^{2,u}(L_{\varphi }^{2})\to l^2(L_{\varphi }^{2})}<\infty.
$$
 It follows from \cite[Theorem 5.1.2]{BL78} that
$$
\left[l^2(L^1_\varphi), l^2(L^2_\varphi)\right]_\theta=l^2\left([L^1_\varphi, L^2_\varphi]_\theta\right).
$$
Set $X_0=X_1=l^2$ and $Y_0= L^1_\varphi, Y_1= L^2_\varphi$. It is well-known that $(l^2)^{\ast}$ is a type 2 space. And $L^p_\varphi$ are complex Banach lattices for $p\ge 1$. Furthermore, for $f_j\in L^p_\varphi$ with $1\le p\le 2$, $j=1, 2, \cdots, m$,   we apply  Minkowski's inequality to get
\beqm
   \left( \sum_{j=1}^m \|f_j\|_{L^p_\varphi}^2\right)^{\frac 12}
&=&  \left( \sum_{j=1}^m \left(\int_{{\mathbb C}^n}\left|f_j e^{-\varphi }\right|^p dv \right)^ {\frac 2 p }\right)^{\frac 12}\\
&\le & \left( \int_{{\mathbb C}^n}\left( \sum_{j=1}^m |f_j|^2\right)^{\frac p 2 }e^{-p\varphi } dv  \right)^{\frac 1p}
= \left\|\left( \sum_{j=1}^m \left|f_j \right|^2\right)^{\frac 1 2 } \right\|_{L^p_\varphi}.
\eqm
This implies that both $L^1_\varphi, L^2_\varphi$ are 2-concave Banach spaces. Applying Proposition \ref{Ko91-1}, we obtain
$$
\left[l^{2,u}(L^1_\varphi), l^{2,u}(L^2_\varphi)\right]_\theta
=\left[l^2\widetilde{\otimes}_\varepsilon L^1_\varphi, l^2\widetilde{\otimes}_\varepsilon L^2_\varphi\right]_\theta=l^2\widetilde{\otimes}_\varepsilon[L^1_\varphi, L^2_\varphi]_\theta.
$$
By interpolation we have
$$
\widehat{T}: l^2\widetilde{\otimes}_\varepsilon[L^1_\varphi, L^2_\varphi]_\theta\to l^2\left([L^1_\varphi, L^2_\varphi]_\theta\right)
$$
is bounded with the norm estimate
$$
\|\widehat{T}\|_{l^2\widetilde{\otimes}_\varepsilon[L^1_\varphi, L^2_\varphi]_\theta\to l^2\left([L^1_\varphi, L^2_\varphi]_\theta\right)}\le  M_1^{1-\theta}M_2^\theta,
$$
where $M_j=\|\widehat{T}\|_{l^2\widetilde{\otimes}_\varepsilon L^j_\varphi\to l^2(L^j_\varphi)},$ $j=1,2.$ Therefore, (\ref{342}) along with (\ref{352}) gives
\beqm
\pi_2(T_\mu)\lesssim \|\widehat{T}\|_{l^2\widetilde{\otimes}_\varepsilon[L^1_\varphi, L^2_\varphi]_\theta\to l^2\left([L^1_\varphi, L^2_\varphi]_\theta\right)}
\le M_1^{1-\theta}M_2^\theta
\lesssim  \|\widetilde{\mu}\|_{L^2},
\eqm
which indicates $T_\mu: F_{\varphi }^{p}\to F_{\varphi }^{p}$ is 2-summing, and  is also $r$-summing  for all $r\ge 1$ with the estimate
\begin{equation}\label{proof-g}
\pi_r(T_\mu)\lesssim \|\widetilde{\mu}\|_{L^2}.
\end{equation}

Now we prove the necessity. Since the cotype of \({F}_{\varphi}^{p}\) is 2 for \(1\le p \le  2\),  we have $T_\mu\in \Pi_r(F_{\varphi}^{p}, {F}_{\varphi }^p)=\Pi_p(F_{\varphi}^{p}, {F}_{\varphi }^p)$. By  Pietsch's domination theorem, there exists a regular probability measure $\sigma$ on ${B}_{{F}_{\varphi }^{p'}}$ such that for each $f\in F_{\varphi}^{p}$, 
 \begin{equation}\label{42}
\left\|T_\mu\left(f\right)\right\|_{{F}_{\varphi}^{p}}^{p}\leq  \big({\pi }_{p}\left(T_\mu \right)\big)^{p}\cdot\int_{  B_{F_\varphi^{p'}}}{\left|\left\langle g, f \right\rangle_{\varphi}\right| }^{p}{d\sigma(g)}.
 \end{equation}
Now for $\delta>0 $  small enough, fix some  $\delta$-lattice $\{z_j\}$.
For $\{c_j\}_{j=1}^m$, set
 $$f_t=\sum_{j=1}^mc_jr_j(t)k_{\varphi,z_j},$$
 where  $r_j(t)$ is the Rademacher function on $[0,1]$.
On one hand, we have $\widetilde{\mu}(z_j)\simeq |T_\mu(k_{\varphi,z_j})(z_j)|e^{-\varphi(z_j)}$. Then Fubini's theorem and Khintchine's inequality give
\beqm
\int_{0}^{1}\left\|T_\mu\left(f_t\right)\right\|_{{F}_{\varphi}^{p}}^{p}dt
&=&\int_{\mathbb{C}^n}\left(\int_{0}^{1}\left|\sum^{m}_{k=1}c_kr_k(t)T_\mu(k_{\varphi,z_k})(z)\right|^{p}dt\right)e^{-p\varphi(z)}dv(z)\\
&\simeq&\int_{\mathbb{C}^n}\left(\sum^{m}_{k=1}\left|c_k\right|^2|T_\mu(k_{\varphi,z_k})(z)|^2\right)^{\fr p 2}e^{-p\varphi(z)}dv(z)\\
&\simeq&\sum^{\infty}_{j=1}\int_{B(z_j,\delta)}\left(\sum^{m}_{k=1}\left|c_k\right|^2|T_\mu(k_{\varphi,z_k})(z)|^2\right)^{\fr p 2}e^{-p\varphi(z)}dv(z)\\
&\gtrsim&\sum^{m}_{j=1}\left|c_j\right|^p\int_{B(z_j, \delta)}|T_\mu(k_{\varphi,z_j})(z)|^pe^{-p\varphi(z)}dv(z)\\
&\gtrsim&\sum^{m}_{j=1}\left|c_j\right|^{p}\widetilde{\mu}(z_j)^p.
\eqm
On the other hand, by using H\"{o}lder's inequality with  $\fr {p'} 2$ and its conjugate exponent, we have
\beqm
&&\int_{0}^{1}\int_{B_{F_\varphi^{p'}}}\left|\left\langle g, f_t \right\rangle_{\varphi} \right|^{p}{d\sigma(g)}dt\\
&=&\int_{B_{F_\varphi^{p'}}}d\sigma(g)\int_{0}^{1}\left|\sum_{j=1}^mc_jr_j(t)\frac{g(z_j)}{\sqrt{K_\varphi(z_j,z_j)}}\right|^{p}dt\\
&\simeq&\int_{B_{F_\varphi^{p'}}}\left(\sum_{j=1}^m|c_j|^2\left|g(z_j)e^{-\varphi(z_j)}\right|^2\right)^{\fr p 2}d\sigma(g)\\
&\lesssim& \left(\sum_{j=1}^m|c_j|^{\frac{2p'}{p'-2}}\right)^{\frac{p'-2}{p'}\cdot \fr p 2}\cdot \int_{B_{F_\varphi^{p'}}}\left(\sum_{j=1}^m\left|g(z_j)e^{-{\varphi(z_j)}}\right|^{p'}\right)^{\frac{p}{p'}}d\sigma(g).
\eqm
For $g\in B_{F_\varphi^{p'}}$, \cite[Proposition 2.3]{SV12} gives
\beqm
\sum_{j=1}^m\left|g(z_j)e^{-{\varphi(z_j)}}\right|^{p'}
&\lesssim& \sum_{j=1}^m\int_{B(z_j,\delta)}\left|g(z)e^{-{\varphi(z)}}\right|^{p'}dv(z)\\
&\lesssim& \int_{\mathbb{C}^n}\left|g(z)e^{-{\varphi(z)}}\right|^{p'}dv(z)\\
&\le & 1.
\eqm
Hence, (\ref{42}) implies that
\begin{equation}\label{m-sum}
\sum^{m}_{j=1}\left|c_j\right|^{p}\widetilde{\mu}(z_j)^p\lesssim \left(\pi_p( T_\mu) \right)^{p}\cdot\left(\sum_{j=1}^m|c_j|^{\frac{2p}{2-p}}\right)^{\frac{2-p}{2}},
\end{equation}
A careful check shows  the constants  $C$ appearing in the estimate ``$\lesssim$'' above are independent of $m$. Therefore, for $\{c_j\}\in l^{\frac{2p}{2-p}}$, there holds
$$
\sum^{\infty}_{j=1}\left|c_j\right|^{p}\widetilde{\mu}(z_j)^p\lesssim \left(\pi_p( T_\mu) \right)^{p}\cdot\|\{|c_j|^p\}\|_{l^{\frac{2}{2-p}}}.
$$
A duality argument with respect to $\frac{2}{2-p}$ and its conjugate $\frac{2}{p}$ indicates
$$
\left[\sum^{\infty}_{j=1}\left(\widetilde{\mu}(z_j)^p\right)^{\frac{2}{p}}\right]^{\frac{p}{2}}\lesssim \left(\pi_p( T_\mu) \right)^{p}.
$$
Therefore,
$$
\left(\sum^{\infty}_{j=1}\widetilde{\mu}(z_j)^{2}\right)^{\frac{1}{2}}\lesssim {\pi }_{p}{\left( T_\mu\right) }.
$$
Notice that $\widehat{\mu}_\delta(z_j)\lesssim \widetilde{\mu}(z_j)$ for small enough $\delta>0$, which together with Lemma \ref{lem2.3} shows $\widetilde{\mu}\in L^ {2}$. Moreover,
\begin{equation}\label{norm-es}
\|\widetilde{\mu}\|_{L^ {2}}\simeq \|\{\widehat{\mu}_\delta(z_j)\}\|_{l^ {2}}\lesssim {\pi }_{p}( T_\mu)\simeq {\pi }_{r}( T_\mu).
\end{equation}

The estimates (\ref{norm-es}) and (\ref{proof-g}) yield the norm equivalence $$
{\pi }_{r}( T_\mu)\simeq \|\widetilde{\mu}\|_{L^ {2}}.
$$
We have obtained  the conclusion  of Theorem \ref{THM2} for $1\le p \le 2$.

\subsection{Case 2: $p\ge 2$ and $1\le r\le p'$.}

  If \({T}_{\mu} : {F}_{\varphi}^{p} \rightarrow  {F}_{\varphi }^{p}\) is \(r\)-summing, then  clearly $T_\mu\in \Pi_{p'}(F_{\varphi}^{p}, {F}_{\varphi }^p)$. Moreover, since $\Pi_{p'}(F_{\varphi}^{p}, {F}_{\varphi }^p)\subset  {\mathcal  B}(F^p_\varphi, F^p_\varphi)$, we have $(T_\mu)^*= T_\mu \in {\mathcal B}(F^{p'}_\varphi, F^{p'}_\varphi)$.  Consider
$$
L^p_\varphi\stackrel{P_\varphi}{\longrightarrow} F^p_\varphi\stackrel{T_\mu}{\longrightarrow} F^p_\varphi,\ \textrm{and} \ \ F^{p'}_\varphi\stackrel{T_\mu}{\longrightarrow} F^{p'}_\varphi  \stackrel{\mathrm{Id}_4} {\longrightarrow} L^{p'}_\varphi.
$$
Set $Q = T_\mu\circ P_\varphi: L_{\varphi }^{p}\to F_{\varphi }^{p}$ and $S = {\rm Id}_4\circ {T}_{\mu}: {F}_{\varphi}^{p'} \rightarrow  L_{\varphi }^{p'}$. It is easy to see that $Q= S^*$ and  $Q\in {\Pi_{p'}} (L_{\varphi}^{p}, F_{\varphi }^{p})$.
 It follows from \cite[Corollary 5.21]{DJT95} that \(S: {F}_{\varphi}^{p'} \rightarrow  L_{\varphi }^{p'}\) is order bounded, i.e., $\mathop{\sup }\limits_{{g \in  {B}_{{F}_{\varphi }^{p'}}}}{\left| T_\mu(g)(z) \right| }\in L_{\varphi}^{p'}$ with the estimate
 \begin{equation}\label{41}
\|\mathop{\sup }\limits_{{g \in  {B}_{{F}_{\varphi }^{p'}}}}{\left| T_\mu(g)(z) \right| }\big\|_{L^{p'}_\varphi}\le \pi_{p'}(T_\mu).
 \end{equation}
Since $ \|k_{\varphi, z}\|_{F^{p'}_\varphi} \simeq 1$, we have
\beqm
\mathop{\sup }\limits_{{g \in  {B}_{{F}_{\varphi }^{p'}}}}{\left| T_\mu(g)(z) \right| } \gtrsim\left| T_\mu(k_{\varphi, z})(z) \right|
 \simeq  \widetilde{\mu }(z)\cdot e^{\varphi(z)}.
\eqm
This together with (\ref{41}) yields
 \begin{equation}\label{4.3}
 \|\widetilde{\mu}\|_{L^{p'}}\lesssim\pi_{p'}(T_\mu)\le \pi_{r}(T_\mu).
 \end{equation}

Conversely, if \(\widetilde{\mu} \in  {L}^{p'}\), Lemma \ref{lem2.3} (2) indicates that \({T}_{\mu} : {F}_{\varphi}^{p} \rightarrow  {F}_{\varphi }^{1}\) is bounded with $\|T_\mu\|\simeq \|\widetilde{\mu}\|_{L^{p'}}$. Meanwhile, since ${\rm Id}_1: {F}_{\varphi}^{1} \rightarrow  {F}_{\varphi }^{2}$ is bounded, as shown in (\ref{1-to-2}) we know ${\rm Id}_1: {F}_{\varphi}^{1} \rightarrow  {F}_{\varphi }^{2}$ is 1-summing with $\pi_1({\rm Id}_1: {F}_{\varphi}^{1} \to {F}_{\varphi }^{2})\lesssim\|{\rm Id}_1\|$. ${\rm Id}_5: {F}_{\varphi}^{2} \rightarrow  {F}_{\varphi }^{p}$ is bounded with $\|{\rm Id}_5\|\le 1$. Consider the commutative diagram:
$$
			\xymatrix{
		 \ \ \ \ 		{F}_{\varphi}^{p} \ \ \ \ 	\ar[r]^{T_\mu} \ar[d]^{T_\mu} &  \ \ \ \ {F}_{\varphi}^{p} \ \ \ \ \\
			 \ \ \ \ 	{F}_{\varphi}^{1} \ \ \ \ 	\ar[r]_{{\rm Id}_1} &  \ \ \ \ F_{\varphi}^{2} \ \ \ \ . \ar[u]_{{\rm Id}_5}}
 $$
Then (\ref{ideal-a}) shows
$
{T}_{\mu}({F}_{\varphi}^{p} \rightarrow  {F}_{\varphi }^{p})={\rm Id}_5 \circ  {\rm Id}_1\circ {T}_{\mu}({F}_{\varphi}^{p} \rightarrow  {F}_{\varphi }^{1})
$
is   $r$-summing with the estimate
 \begin{equation}\label{4.4}
\pi_r(T_\mu)\le\pi_1(T_\mu)\lesssim \|\widetilde{\mu}\|_{L^{p'}}.
 \end{equation}

 Combining (\ref{4.3}) and (\ref{4.4}), we obtain the desired norm equivalence in (\ref{main-thm}) for $p\ge 2$ and $1\le r\le p'$.

\subsection{Case 3: $p\ge 2$ and $p' \leq  r \leq p$.}

Suppose \(\widetilde{\mu} \in  {L}^{r}\). Consider
$$ F^\infty_\varphi\stackrel{T_\mu}{\longrightarrow} F^r_\varphi\stackrel{\mathrm{Id}_6} {\longrightarrow} L^{r}_\varphi\stackrel{P_\varphi}{\longrightarrow} F^{r}_\varphi.
$$
We claim \(\mathrm{Id}_6\circ{T}_{\mu} : {F}_{\varphi}^{\infty} \rightarrow  L_{\varphi }^{r}\) is order bounded. In fact, \cite[Lemma 2.2]{HL14} and H\"{o}lder's inequality give
\beqm
\mathop{\sup }\limits_{{g \in  {B}_{{F}_{\varphi }^{\infty}}}}{\left| T_\mu(g)(z) \right| }
&\lesssim& \mathop{\sup }\limits_{{g \in  {B}_{{F}_{\varphi }^{\infty}}}}\int_{\mathbb{C}^n}|g(\xi)||K_\varphi(z,\xi)|e^{-2\varphi(\xi)}\widehat{\mu}_\delta(\xi)dv(\xi)\\
&\le & \int_{\mathbb{C}^n}|K_{\varphi,z}(\xi)|e^{-\varphi(\xi)}\widehat{\mu}_\delta(\xi)dv(\xi)\\
&\le & \left(\int_{\mathbb{C}^n}|K_{\varphi,z}(\xi)|e^{-\varphi(\xi)}\widehat{\mu}_\delta(\xi)^rdv(\xi)\right)^{\fr 1 r}\cdot \|K_{\varphi,z}\|_{F_\varphi^1}^{\fr 1 {r'}}.\\
\eqm
Hence
\beqm
&&\left(e^{-\varphi(z)}\cdot\mathop{\sup }\limits_{{g \in  {B}_{{F}_{\varphi }^{\infty}}}}{\left| T_\mu(g)(z) \right|}\right)^r\\
&\lesssim&
 \left(e^{-\varphi(z)}\|K_{\varphi,z}\|_{F_\varphi^1}\right)^{\fr r {r'}}
\left(e^{-\varphi(z)}\int_{\mathbb{C}^n}|K_{\varphi,z}(\xi)|e^{-\varphi(\xi)}\widehat{\mu}_\delta(\xi)^rdv(\xi)\right)\\
&\lesssim&
e^{-\varphi(z)}\int_{\mathbb{C}^n}|K_{\varphi,z}(\xi)|e^{-\varphi(\xi)}\widehat{\mu}_\delta(\xi)^rdv(\xi).
\eqm
Integrating both sides over $\mathbb{C}^n$ against the measure $dv$, applying Fubini's theorem,  we have
\beqm
&& \int_{\mathbb{C}^n}\left(e^{-\varphi(z)}\cdot\mathop{\sup }\limits_{{g \in  {B}_{{F}_{\varphi }^{\infty}}}}{\left| T_\mu(g)(z) \right|}\right)^rdv(z)\\
&\lesssim&
\int_{\mathbb{C}^n}\widehat{\mu}_\delta(\xi)^rdv(\xi)
\int_{\mathbb{C}^n}e^{-\varphi(z)-\varphi(\xi)}|K_{\varphi}(z,\xi)|dv(z)\\
&\lesssim& \int_{\mathbb{C}^n}\widehat{\mu}_\delta(\xi)^rdv(\xi).
\eqm
This shows \({\textrm{Id}_6} \circ T_{\mu} : {F}_{\varphi}^{\infty} \rightarrow  {L}_{\varphi }^{r}\) is order bounded, and hence it is $r$-summing with
$$
\pi_r({\textrm{Id}_6} \circ T_\mu)\le \big\|\mathop{\sup }\limits_{{g \in  {B}_{{F}_{\varphi }^{\infty}}}}{\left| T_\mu(g) \right| }\big\|_{L^{r}_\varphi}\lesssim \|\widehat{\mu}_\delta\|_{L^{r}}.
$$
Thus we get
$
T_\mu({F}_{\varphi}^{\infty} \to F_{\varphi }^{r})=P_\varphi(L_{\varphi}^r \to F_{\varphi }^{r})\circ\left( {\textrm{Id}_6} \circ T_\mu\right)
$
is $r$-summing with
$$
\pi_r(T_\mu: {F}_{\varphi}^{\infty} \to F_{\varphi }^{r})\lesssim \|\widehat{\mu}_\delta\|_{L^{r}}.
$$
Consider the commutative diagram:
$$
			\xymatrix{
				 \ \ \ \ {F}_{\varphi}^{p} \ \ \ \ 	\ar[r]^{T_\mu} \ar[d]^{{\rm Id}_7} & \ \ \ \ {F}_{\varphi}^{p} \ \ \ \ \\
				 \ \ \ \ {F}_{\varphi}^{\infty} \ \ \ \ 	\ar[r]_{T_\mu} & \ \ \ \ F_{\varphi}^{r} \ \ \ \ . \ar[u]_{{\rm Id}_8}}
 $$
Notice that both ${\rm Id}_7: {F}_{\varphi}^{p} \rightarrow  {F}_{\varphi }^{\infty}$ and ${\rm Id}_8: {F}_{\varphi}^{r} \rightarrow  {F}_{\varphi }^{p}$  are bounded, then (\ref{ideal-a}) gives \({T}_{\mu} : {F}_{\varphi}^{p} \rightarrow  {F}_{\varphi }^{p}\) is $r$-summing with
 \begin{equation}\label{4.22}
\pi_r(T_\mu : {F}_{\varphi}^{p} \rightarrow  {F}_{\varphi }^{p})\lesssim \|\widehat{\mu}_\delta\|_{L^{r}}\simeq \|\widetilde{\mu}\|_{L^{r}}.
 \end{equation}

Conversely, if \({T}_{\mu} : {F}_{\varphi}^{p} \rightarrow  {F}_{\varphi }^{p}\) is $r$-summing, set
$$
F: {\mathbb{C}}^{n}  \longrightarrow  {F}_{\varphi}^{p},
$$
$$
\ \ \ \ \ \ \ \ \ \ z   \longmapsto  k_{\varphi,z}(\cdot)
$$
From Lemma \ref{lem8} we know, for the above $F$ on measurable space $(\mathbb{C}^{n}, \upsilon)$, there holds
$$
\int_{\mathbb{C}^{n} }\|T_\mu  (k_{\varphi,z})\|_{{F}_{\varphi}^{p}}^{r}d\upsilon(z) \leq  \left({\pi }_{r}{\left( T_\mu\right) }\right)^{r}\cdot\mathop{\sup }\limits_{{g  \in  {B}_{{F}_{\varphi}^{p'}}}}\int_{\mathbb{C}^{n}}{\left|\langle g, k_{\varphi,z} \rangle_\varphi\right| }^{r}{dv(z)}.
$$
On one hand, the sub-mean value property of $|T_\mu (k_{\varphi,z})e^{-\varphi(\cdot)}|^p$ yields
\beqm
\|T_\mu (k_{\varphi,z})\|_{{F}_{\varphi}^{p}}^{p}
&\ge & \int_{B(z,s)}|T_\mu (k_{\varphi,z})(\xi)|^pe^{-p\varphi(\xi)}dv(\xi)\\
&\gtrsim& |T_\mu (k_{\varphi,z})(z)|^pe^{-p\varphi(z)}\\
&\simeq& \widetilde{\mu}(z)^p.
\eqm
On the other hand, since $r\ge p'$, we have
\beqm
\mathop{\sup }\limits_{{g  \in  {B}_{{F}_{\varphi}^{p'}}}}\int_{\mathbb{C}^{n}}{\left|\langle g, k_{\varphi,z} \rangle_\varphi\right| }^{r}{dv(z)}
&=& \mathop{\sup }\limits_{{g  \in  {B}_{{F}_{\varphi}^{p'}}}}\int_{\mathbb{C}^{n}}{\left|g(z)e^{-\varphi(z)}\right| }^{r}{dv(z)} \\
&\lesssim &  \mathop{\sup }\limits_{{g  \in  {B}_{{F}_{\varphi}^{p'}}}}\|g\|_{F_\varphi^{p'}}^{r}=1.
\eqm
The above estimates show that
$$
\int_{\mathbb{C}^{n}}\widetilde{\mu}(z)^{r}{dv(z)}\lesssim \left({\pi }_{r}{( T_\mu)}\right)^{r},
$$
which tells us \(\widetilde{\mu}\in  {L}^{r}\) with $\|\widetilde{\mu}\|_{L^r}\lesssim \pi_r(T_\mu)$.
This together with (\ref{4.22}) gives the  norm equivalence we want.

\subsection{Case 4: $p\ge 2$ and $r \geq p$.}

If \(\widetilde{\mu} \in  {L}^{p}\), then the conclusion of Case 3 in (\ref{main-thm}) shows that \({T}_{\mu} : {F}_{\varphi}^{p} \rightarrow  {F}_{\varphi }^{p}\) is \(p\)-summing. Consequently, for any $r \geq p$,   \({T}_{\mu} : {F}_{\varphi}^{p} \rightarrow  {F}_{\varphi }^{p}\) is \(r\)-summing with 
\begin{equation}\label{case-4-1}
\pi_r(T_\mu )\le \pi_p(T_\mu )\simeq \|\widetilde{\mu}\|_{L^{p}}.
\end{equation}

Conversely, suppose \({T}_{\mu} : {F}_{\varphi}^{p} \rightarrow  {F}_{\varphi }^{p}\) is \(r\)-summing. Since the cotype of ${F}_{\varphi}^{p}$ is $p \geq  2$, it follows from \cite[Theorem 11.13]{DJT95} that $T_\mu: F_\varphi^p\to  F_{\varphi}^p$ is $(p, 2)$-summing (see \cite[p.197]{DJT95} for the definition of $(r,s)$-summing operators). Moreover, a careful check of the proof of \cite[Theorem 11.13]{DJT95}  shows that 
\beqm
\left(\sum_{j=1}^{m}\|T_\mu f_j\|_{{F}_{\varphi}^{p}}^{p}\right)^{\fr 1 p}
&\lesssim& \pi_r(T_\mu) \cdot \mathop{\sup }\limits_{{g  \in  {B}_{{F}_{\varphi}^{p'}}}}\left(\sum_{j=1}^{m}\left|\langle g, f_j \rangle_\varphi\right|^2\right)^{\fr 1 2}\\
\eqm
for $\{f_j\}_{j=1}^m\subset F_{\varphi}^{p}$. Since $p'\le 2$, with the elementary inequality that $(a+b)^\alpha\le a^\alpha + b^\alpha$ for $a, b\ge 0$ and $0<\alpha \le 1$ we know
\beqm
\left(\sum_{j=1}^{m}\|T_\mu f_j\|_{{F}_{\varphi}^{p}}^{p}\right)^{\fr 1 p}
&\lesssim &  \pi_r(T_\mu)\cdot \mathop{\sup }\limits_{{g  \in  {B}_{{F}_{\varphi}^{p'}}}}\left(\sum_{j=1}^{m}\left|\langle g, f_j \rangle_\varphi\right|^{p'}\right)^{\fr 1 {p'}}.
\eqm
Now take some  $\delta$-lattice   $\{z_j\}$ with $\delta>0$ small enough, set $f_j=k_{\varphi, z_j}$. Then  we obtain
\beqm
\left(\sum_{j=1}^{m}\widetilde{\mu}(z_j)^{p}\right)^{\fr 1 p}&\lesssim&\left(\sum_{j=1}^{m}\|T_\mu (k_{\varphi, z_j})\|_{{F}_{\varphi}^{p}}^{p}\right)^{\fr 1 p}\\
&\lesssim& \pi_r(T_\mu)\cdot \mathop{\sup }\limits_{{g  \in  {B}_{{F}_{\varphi}^{p'}}}}\left(\sum_{j=1}^{m}\left|\langle g, k_{\varphi, z_j} \rangle_\varphi\right|^{p'}\right)^{\fr 1 {p'}}\\
&\simeq&\pi_r(T_\mu)\cdot \mathop{\sup }\limits_{{g  \in  {B}_{{F}_{\varphi}^{p'}}}}\left(\sum_{j=1}^{m}\left|g(z_j)e^{-\varphi(z_j)} \right|^{p'}\right)^{\fr 1 {p'}}\\
&\lesssim& \pi_r(T_\mu)\cdot \mathop{\sup }\limits_{{g  \in  {B}_{{F}_{\varphi}^{p'}}}}\|g\|_{{F}_{\varphi}^{p'}}^{p'}\\
&= & \pi_r(T_\mu).
\eqm
Similar to (\ref{norm-es}) we get \(\widetilde{\mu} \in  {L}^{p}\) and
\begin{equation}\label{case-4-2}
\|\widetilde{\mu}\|_{L^{p}}\lesssim \pi_r(T_\mu).
\end{equation}

Combining (\ref{case-4-1}) and (\ref{case-4-2}) we have (\ref{main-thm}) for $p\ge 2$ and $r \geq p$. The proof is finished.

\section{Further Remarks}

It would be interesting  to extend our investigation to the mapping properties of Toeplitz operators between different Fock spaces. More precisely, given
$1\le p, q <\infty$, one may ask for necessary and sufficient conditions on the measure $\mu$  that guarantee the operator
$$
   T_\mu: F^p_\varphi \to F^q_\varphi
$$
is $r$-summing. Such results would generalize the self-mapping case treated here, and could lead to a more flexible operator-theoretic framework encompassing interpolation phenomena and cross-space summability.

Let $\mu$ be a positive Borel measure on \(\mathbb{C}^n\). As an easy consequence we know $T_\mu \in \Pi_1(F^1_\varphi, F^2_\varphi)$ if and only if $T_\mu \in \mathcal{B} (F^1_\varphi, F^2_\varphi)$. Then from \cite{HL14} we see that $T_\mu \in \Pi_1(F^1_\varphi, F^2_\varphi)$ if and only if
$ \widetilde{\mu} \in L^\infty$ and ${\pi }_{r}{\left( T_\mu: {F}_{\varphi}^{1} \rightarrow  {F}_{\varphi }^{2}\right) }\simeq \|\widetilde{\mu}\|_{L^ {\infty}}$, $1\le r< \infty$. For $1<p<2$ we have the following proposition.

 \begin{proposition} For  \(1< p < 2\) and $1\le r<\infty$, then  \({T}_{\mu} \in \Pi_r({F}_{\varphi}^{p}, {F}_{\varphi }^{2})\) if and only if \(\widetilde{\mu} \in L^ {\frac{p}{p-1}}\). Furthermore,
 $$
{\pi }_{r}{\left( T_\mu: {F}_{\varphi}^{p} \rightarrow  {F}_{\varphi }^{2}\right) }\simeq \|\widetilde{\mu}\|_{L^ {\frac{p}{p-1}}}.
 $$
\end{proposition}

\begin{proof}
Suppose \(\widetilde{\mu} \in L^ {\frac{p}{p-1}}\) for \(1< p < 2\), it follows from Lemma \ref{lem2.3} (2) that ${T}_\mu: {F}_{\varphi}^{p} \rightarrow  {F}_{\varphi }^{1}$ is bounded with $\|T_\mu\|_{F_\varphi^p \rightarrow  {F}_{\varphi }^{1}}\lesssim \|\widetilde{\mu}\|_{L^ {\frac{p}{p-1}}}$. Meanwhile,  ${\rm Id}: {F}_{\varphi}^{1} \rightarrow  {F}_{\varphi }^{2}$ is 1-summing with $\pi_1({\rm Id}: F_\varphi^1 \rightarrow  F_\varphi ^2)\lesssim \|{\rm Id}\|$. Therefore,
$$
{T}_{\mu}({F}_{\varphi}^{p} \rightarrow  {F}_{\varphi }^{2})={\rm Id}({F}_{\varphi}^{1} \rightarrow  {F}_{\varphi }^{2})\circ T_\mu(F_\varphi^p \rightarrow  {F}_{\varphi }^{1})
$$
is \(1\)-summing with
$$
\pi_1({T}_{\mu}:{F}_{\varphi}^{p} \rightarrow  {F}_{\varphi }^{2})\lesssim\pi_1({I}_{d}:{F}_{\varphi}^{1} \rightarrow  {F}_{\varphi }^{2})\cdot\|T_\mu\|_{F_\varphi^p \rightarrow  {F}_{\varphi }^{1}}\lesssim\|\widetilde{\mu}\|_{L^ {\frac{p}{p-1}}}.
$$

Conversely, if \({T}_{\mu} : {F}_{\varphi}^{p} \rightarrow  {F}_{\varphi }^{2}\) is \(r\)-summing, with the same  approach  as in the proof of the necessity for Case 1, we know  \(\widetilde{\mu} \in L^ {\frac{p}{p-1}}\) with  $
\|\widetilde{\mu}\|_{L^ {\frac{p}{p-1}}}\lesssim {\pi }_{r}{\left( T_\mu: {F}_{\varphi}^{p} \rightarrow  {F}_{\varphi }^{2}\right) }.$
\end{proof}

Our investigation of
$r$-summing positive Toeplitz operators on Fock spaces also points toward potential progress in the study of Bergman spaces on the unit ball. Since both frameworks rely fundamentally on reproducing kernel estimates and Carleson measure techniques, it is natural to expect that the characterization of absolutely summing operators on
$F^p_\varphi$ may serve as a model for analogous problems in
$A^p$. In particular, the structural insights obtained here could help identify precise measure-theoretic conditions guaranteeing the absolute summability of
$T_\mu $
 acting on Bergman spaces, thereby opening a parallel line of investigation in the ball setting. Although one could, in principle, formulate a conjectural description for the full range
$1<p<\infty$, the general case leads to considerably more intricate expressions. For the sake of clarity, we confine ourselves here to the range
$1<p\le 2$ and state the following conjecture.

\begin{conjecture}
Suppose $\mu$ is a finite positive Borel measure on the unit ball $\mathbb B$ of ${\mathbb C}^n$ and $1<p\le 2$. Then the Toeplitz operator $T_\mu$ is $r$-summing on $A^p(\mathbb B)$ if and only if
$$
      \int _{\mathbb B} \widetilde{\mu}(z)^2K(z, z) dv(z)<\infty,
$$
where $\widetilde{\mu}$ is the Berezin transform of $\mu$ and $K(\cdot, \cdot)$ is the Bergman kernel on $A^2(\mathbb B)$.
\end{conjecture}

\textbf{Acknowledgement:} The authors are grateful to the referee for  valuable comments and helpful suggestions. Hu is supported by the National Natural Science
Foundation of China under Grant No.  12071130, 12171150. Wang is supported by the National Natural Science Foundation of China under Grant No. 12671152 and Yanling Youqing Project of Lingnan Normal University.

\begin{center}

\noindent
\end{center}


\begin{thebibliography}{99}

\bibitem{BC86}  C. Berger, L. Coburn, Toeplitz operators and quantum mechanics, J. Funct. Anal., 68, 273-299  (1986).

\bibitem{BC87}  C. Berger, L. Coburn, Toeplitz operators on the Segal-Bargmann space, Trans. Amer. Math. Soc., 301, 813-829  (1987).

\bibitem{BL78} J. Bergh and J. L\"{o}fstr\"{o}m, Interpolation spaces, Springer-Verlag, 1978.


\bibitem{CHW24} J. Chen, B. He,  M. Wang, Absolutely summing Carleson embeddings on weighted Fock spaces with $A_{\infty}$-type weights, J. Oper. Theory, 95, 119-151 (2026).



\bibitem{DC00} A. Defant, M. Carsten, A complex interpolation formula for tensor products of vector-valued Banach function spaces. Arch. Math., 74,  441-451 (2000).


\bibitem{DJT95} J. Diestel, H. Jarchow, A. Tonge, Absolutely summing operators, Cambridge University Press, Cambridge (1995).

\bibitem{FL22} T. Far\`{e}s, P. Lef\`{e}vre, Absolutely summing weighted composition operators on Bloch spaces, J. Oper. Theory, 88, 407-443 (2022).


\bibitem{Gr53} A. Grothendieck, R\'{e}sum\'{e} de la th\'{e}orie m\'{e}trique des produits tensoriels topologiques, Bol. Soc. Mat. S\~{a}o Paulo, 8, 1-79 (1953).

\bibitem{Gr55} A. Grothendieck, Produits tensoriels topologiques et espaces nucl\'{e}aires. Mem. Amer. Math. Soc., 16, 140 (1955).


\bibitem{HJL24} B. He, J. Jreis, P. Lef\`{e}vre, Z. Lou, Absolutely summing Carleson embeddings on Bergman spaces, Adv. Math., 439, 109495 (2024).

\bibitem{HL11} Z. Hu, X. Lv,  Toeplitz operators from one Fock space to another, Integ. Equ.  Oper. Theory, 70, 541-559 (2011).

\bibitem{HL14} Z. Hu, X. Lv, Toeplitz operators on Fock spaces $F^p(\varphi)$, Integr. Equ. Oper. Theory, 80, 33-59 (2014).

\bibitem{HL17} Z. Hu, X. Lv, Positive Toeplitz operators between different doubling Fock spaces, Taiwan. J. Math., 21, 467-487 (2017).


\bibitem{IZ10} J. Isralowitz, K. Zhu, Toeplitz Operators on the Fock Space, Integr. Equ. Oper. Theory, 66, 593-611 (2010).


\bibitem{IVW15} J. Isralowitz, J. Virtanen, L. Wolf, Schatten class Toeplitz operators on generalized Fock spaces, J. Math. Anal. Appl. 421, 329-337 (2015).

\bibitem{JL24} J. Jreis, P. Lef\`{e}vre, Some operator ideal properties of Volterra operators on Bergman and Bloch spaces, Integ. Equa. Oper. Theory, 96, No. 1 (2024).

\bibitem{JMP08} H. Junek, M. Matos,  D. Pellegrino, Inclusion theorems for absolutely summing holomorphic mappings, Proc. Amer. Math. Soc., 136, 3983-3991 (2008).

\bibitem{Ko91} O. Kouba,  On the interpolation of injective or projective tensor products of Banach spaces, J. Funct. Anal., 96, 38-61 (1991).


\bibitem{LR18} P. Lef\`{e}vre, L. Rodr\'{i}guez-Piazza, Absolutely summing Carleson embeddings on Hardy spaces, Adv. Math. 340, 528-587 (2018).

\bibitem{LT79}  J. Lindenstrauss, L. Tzafriri, Classical Banach spaces II: Function spaces, Springer-Verlag, (1979).


\bibitem{OP16} R. Oliver, D. Pascuas, Toeplitz operators on doubling Fock spaces, J. Math. Anal. Appl., 435, 1426-1457 (2016).

\bibitem{Pe77} A. Peeczynski, Banach spaces of analytic functions and absolutely summing operators, American Mathematical Society, Providence, R.I. (1977).

\bibitem{Pi67} A. Pietsch, Absolut $p$-summierende Abbildungen in normierten R\"{a}umen, Studia Math. 28, 333-353 (1967).

\bibitem{Pi72} A. Pietsch, Nuclear locally convex spaces, Springer-Verlag (1972).

 \bibitem{Pi78} A. Pietsch, Operator ideals. VEB Deustscher Verlag der Wissensch., 1978; North-Holland (1980).



\bibitem{Ry02} A. Ryan, Introduction to tensor products of Banach spaces, Springer (2002).

\bibitem{SV12} A. Schuster, D. Varolin, Toeplitz operators and Carleson measures on generalized Bargmann-Fock spaces, Integr. Equ. Oper. Theory, 72, 363-392 (2012).

\bibitem{Wo91} P. Wojtaszczyk, Banach spaces for analysts, Cambridge University Press, Cambridge (1991).

\bibitem{Zh07} K. Zhu, Operator theory in function spaces, American Mathematical Society, New York (2007).

\bibitem{Zh12} K. Zhu, Analysis on Fock spaces, Springer, New York (2012).


\end{thebibliography}
\end{document}